\documentclass[11pt]{article}
\usepackage{graphicx}
\usepackage{amsmath}
\usepackage{amssymb}
\usepackage{amscd}
\usepackage{amsthm}
\usepackage{hyperref}
\usepackage{color}
\usepackage{enumerate}
\usepackage{caption}
\usepackage{subcaption}
\usepackage{natbib}
\usepackage{authblk}
\usepackage[top=1in, bottom=1in, left=1in, right=1in]{geometry}

\DeclareMathOperator{\sech}{sech}

\newtheorem{remark}{Remark}
\newcommand\eJG[1]{{\color{black}{#1}}}
\newcommand\eCJ[1]{{\color{black}{#1}}}
\newcommand\eAR[1]{{\color{black}{#1}}}

\title{Mixed mode oscillations of the El Ni\~no-Southern Oscillation}

\author[1]{Andrew Roberts}
\author[2]{Esther Widiasih}
\author[3]{Axel Timmermann}
\author[4]{Chris K. R. T. Jones}
\author[1]{John Guckenheimer}

\affil[1]{Cornell University, Department of Mathematics}

\affil[2]{University of Hawai'i - West O'ahu, Division of Mathematics and Science}

\affil[3]{University of Hawai'i at Manoa, International Pacific Research Center and Department of Oceanography}

\affil[4]{University of North Carolina, Department of Mathematics}

\begin{document}
\maketitle

\begin{abstract} 
Very strong El Ni\~no events occur sporadically every 10-20 years. The origin of this {\it bursting} behavior still remains elusive. Using a simplified 3-dimensional dynamical model of the tropical Pacific climate system, which captures the El Ni\~no-Southern Oscillation (ENSO) combined with recently developed mathematical tools for fast-slow systems we show that decadal ENSO bursting behavior can be explained as a Mixed Mode Oscillation (MMO), which also predicts a critical threshold
for rapid amplitude growth. It is hypothesized that the MMO dynamics of the low-dimensional climate model can be linked to
a saddle-focus \eJG{equilibrium point},
which mimics a tropical Pacific Ocean state without ocean circulation.
\end{abstract}


\section{Introduction}

The ENSO phenomenon is the dominant source of interannual climate variability. El Ni\~no events are characterized by positive sea surface temperature anomalies (SSTA) in the eastern equatorial Pacific which cause anomalous diabatic heating of the atmosphere. This in turn drives global atmospheric planetary wave adjustments impacting weather in regions far away from the tropical Pacific \citep{bjerknes69, DeserWallace1990, karoly89}.

\begin{figure}[t]
        \centering
        \includegraphics[width=\textwidth]{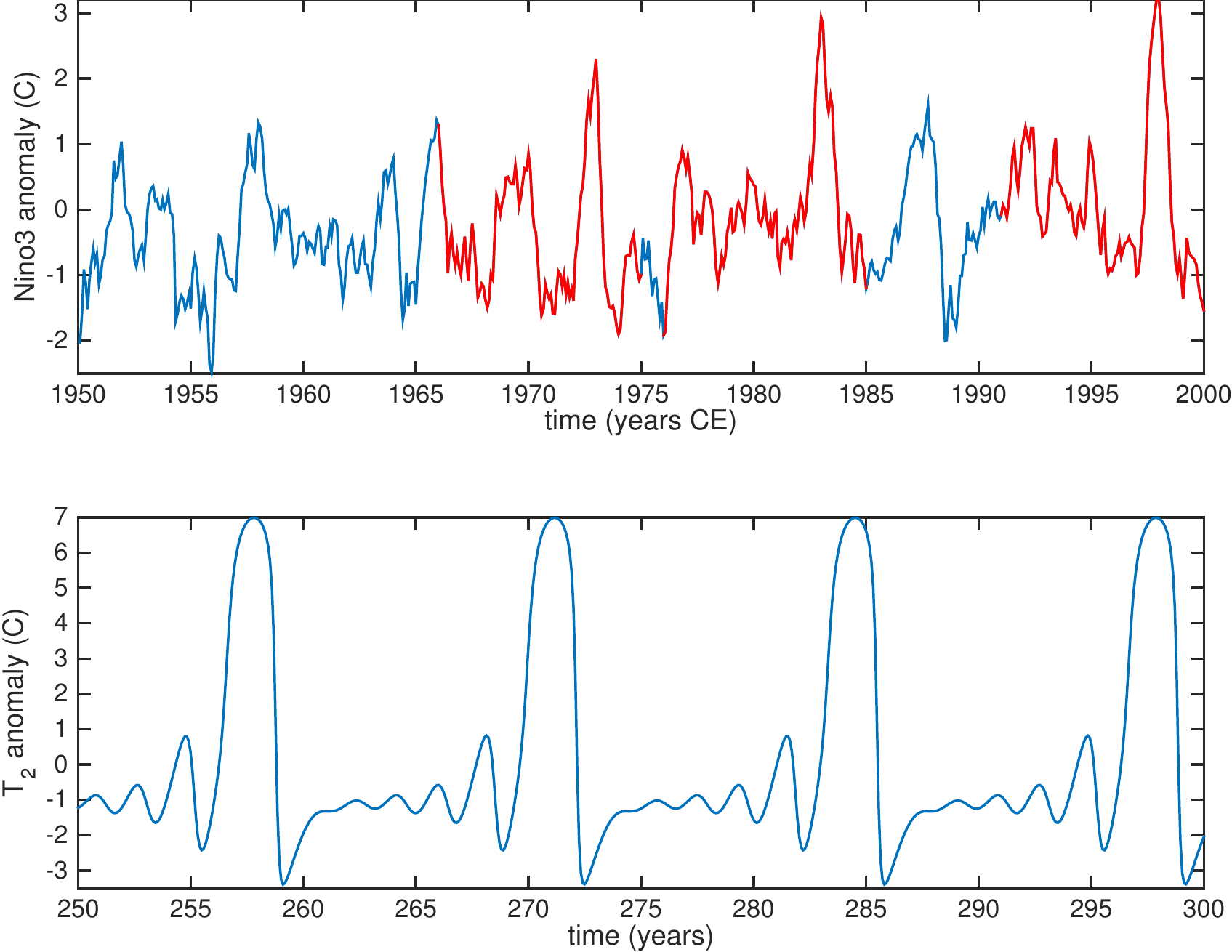}
        \caption{Upper panel: Observed sea surface temperature anomalies averaged over Ni\~no 3 region. High temperatures correspond to strong El Ni\~{n}o events; Red lines indicate a repeated El Ni\~{n}o bursting pattern. 
        Lower panel: Simulated eastern tropical Pacific SSTA in build-up phase and termination phase of strong El Ni\~no events using equations \ref{fullT} and the parameters in Table 1.     
      }
        \label{fig:bursting}
\end{figure}

ENSO can be regarded as a coupled atmosphere/ocean instability, subject to atmospheric noise, external periodic forcing through the annual cycle and atmosphere and ocean nonlinearities. From a dynamical systems perspective the oscillatory ENSO mode emerges from a Hopf bifurcation. To capture the irregular nature of ENSO, its occasional amplitude bursting \citep{tim03} (Figure \ref{fig:bursting}) and the presence of extreme El Ni\~no events, nonlinear extensions to the standard ENSO recharge oscillator theory \citep{jin97,jin98} have been suggested, including noise-induced instabilities \citep{levine2010}, frequency entrainment, chaotic dynamics originating from annual cycle/ENSO interactions \citep{Chang1995, Chang1994, jin94, liu02, stein2014, timmermann02, timmermann03solo, tziperman94} and homoclinic/heteroclinic dynamics \citep{tim03}. 

In view of the socio-economic impacts of very extreme El Ni\~no events, it is paramount to further refine existing theories for the decadal-scale emergence of these large events. Here we hypothesize that ENSO's bursting behavior can be interpreted in terms of Mixed-Mode Oscillations. 
To demonstrate the applicability of this powerful mathematical framework for ENSO, we adopt the same 3-dimensional nonlinear ENSO recharge oscillator equations used \eAR{by \cite{tim03}. }

\eCJ{Mixed-Mode Oscillations (MMOs)} are
 \eCJ{comprised} of $L_1$ large amplitude oscillations (LAOs) followed by $s_1$ small amplitude oscillations (SAOs), then $L_2$ large oscillations, $s_2$ small ones, and so on.  The sequence $L_1^{s_1} L_2^{s_2}L_3^{s_3} \ldots$ is called the MMO signature \citep{mmoSurvey}.  There are a number of geometric explanations for MMOs, and the Shilnikov mechanism mentioned \eAR{by \cite{tim03}} is the poster-child for MMOs in a system 
\eJG{with a single time-scale.}
In systems with multiple time scales, MMO mechanisms can be more robust, and these mechanisms are surveyed by \eAR{\cite{mmoSurvey}}.  The bursting behavior observed in Figure \ref{fig:bursting} is often an indicator that the underlying model has a multiple time-scale structure.  

Here we set out to analyze a dimensionless version of the model from \eAR{\cite{tim03}} to show that there is a parameter regime in which the model exhibits MMOs. \eCJ{We isolate two distinct types of MMO: one has an SAO which is monotonically increasing in amplitude while the other has SAOs that first decrease in amplitude and then increase again. They are each formed by a distinct mathematical mechanism that will be explained in the analysis below. We shall refer to them as monotone and non-monotone SAOs respectively.} We also relate the observed MMOs to the  Shilnikov mechanism discussed \eAR{by \cite{tim03}.}  

\eCJ{We stress that the model under consideration is idealized and omits many processes certainly relevant to the dynamics of ENSO. The goal of the paper is to show that the basic mechanisms captured by the recharge oscillator can account for the bursting seen in ENSO.}

The outline of the paper is as follows: In Section 2 we describe the physical model.  Then we prepare the model for MMO analysis through a coordinate transformation and nondimensionalization.  In particular, we identify the dimensionless model as a multiple time-scale problem and identify the time-scale parameters.  In Section 3 we analyze the model with a focus on demonstrating the existence of MMOs.  This analysis includes an in-depth discussion of the mathematical theory behind both mechanisms leading to MMOs in the ENSO model.  Finally, we conclude with a discussion of the physical implications in Section 4.

\section{A Nonlinear Model for ENSO}
\subsection{The Physical Model}
The fundamental dynamics of ENSO can be described in terms of the recharge oscillator paradigm \citep{jin97}. The key  regions for ENSO physics are the western and eastern equatorial Pacific.  The important dynamical variables are the temperatures in these regions ($T_1,T_2$), eastern tropical Pacific subsurface temperature ($T_{sub}$) and its linkage to western tropical Pacific thermocline changes ($h_1$). Following the original model \citep{jin98,tim03} the underlying ordinary differential equations describing ENSO and its linkage to the mean steady state can be written as  
\begin{equation}  
\label{fullT}
\begin{cases}
	\displaystyle{ \frac{dT_1}{dt} } = -\alpha (T_1-T_r) - \epsilon \mu (T_2-T_1)^2 \\[.6em]
	\displaystyle{ \frac{dT_2}{dt} }= -\alpha (T_2-T_r) + \zeta \mu (T_2-T_1) (T_2-T_{sub}(T_1, T_2, h_1)) \\[.6em]
	\displaystyle{ \frac{dh_1}{dt} }=  r\left(-h_1- \displaystyle{ \frac{bL \mu (T_2-T_1)}{2 \beta} }\right),
\end{cases}
\end{equation}
 where $T_1,$ $T_2$, and $h_1$ represent the equatorial temperature of the Western Pacific, equatorial temperature of the Eastern Pacific, and thermocline depth of the Western Pacific, respectively.  The first terms of the temperature equations represent a relaxation of the system back towards a climatological mean state $T_r$, which represents the radiative-convective equilibrium state (the temperature which the tropical Pacific would attain in the absence of ocean dynamics). The second quadratic term in the $T_1$ equation captures the anomaly wind-driven zonal advection of temperature. The wind anomalies themselves are determined by the east/west temperature gradient ($T_2-T_1$), thus leading to a quadratic dependence in temperature.

 The eastern tropical Pacific temperature tendency is determined by relaxation towards $T_r$ and the vertical temperature advection term, which depends on $T_2$ and the subsurface temperature $T_{sub}$.  The prognostic equation for the western tropical Pacific thermocline variation, $h_1$, captures a relaxation back to mean climatological conditions with a damping timescale $r$ and the effect of wind-stress curl changes on the Sverdrup transport. The latter describes the corresponding discharging (recharging) process for El Ni\~no (La Ni\~na) events, which is assumed to be proportional to the east/west temperature gradient. The system of ordinary differential equations is closed by a nonlinear parameterization of eastern tropical Pacific subsurface temperature variations and \eAR{eastern thermocline depth $h_2=h_1+bL\mu(T_2-T_1)/\beta$: 
\begin{equation*}
	T_{sub}(T_1, T_2, h_1)=\frac{T_r+T_{r_0}}{2} - \frac{T_r-T_{r_0}}{2} \tanh \left( \frac{ H-z_0+h_1+bL\mu(T_2-T_1)/\beta}{h^*}, \right). 
\end{equation*}
The difference $T_r-T_{r_0}$ controls the nonlinear scaling between thermocline anomalies and subsurface temperature anomalies relative to a mean subsurface temperature of $(T_r+T_{r_0})/2$. $T_{r_0}$ corresponds to a mean eastern equatorial temperature attained at a depth of about 75m, $H$ denotes an eastern thermocline reference depth, and $h^*$ indicates the sharpness of the thermocline.}

\begin{table} [h!]
\begin{center}
\addtolength{\tabcolsep}{1mm}
\renewcommand{\arraystretch}{1.2}
\begin{tabular}{|c|c|c||c|c|c|}
\hline
\textbf{Parameter}  & \textbf{Value} & \textbf{Unit} & \textbf{Parameter}  & \textbf{Value} & \textbf{Unit} \\
\hline
\hline 
T$_{r0}$ & 16 & $^\circ$C  & T$_r$ & 29.5 & $^\circ$C   \\
\hline 
$\alpha$ & 1/180 & day$^{-1}$ & r & 1/400 & day$^{-1}$ \\
\hline
H & 100 & meter (m) & z$_0$ & 75 & m \\
\hline
h$^*$ & 62 & m & $\mu$ & 0.0026 & K$^{-1}$ day$^{-1}$ \\
\hline
$\epsilon$& 0.11  &1& $\zeta$ & 1.3 &1 \\
\hline
L & $15 \times 10^6$ & m & $\frac{\mu b L}{\beta}$  & 22 & {mK$^{-1}$} \\
\hline
\end{tabular}
\caption{Parameters used \eAR{by \cite{tim03}}.} \label{physparameter}
\end{center}
\end{table}

The parameters used \eAR{by \cite{tim03}} are listed in Table \ref{physparameter}. To illustrate the dynamical behaviour of this simplified ENSO system we conduct a numerical simulation of the system \eqref{fullT}. The results are shown in the lower panel of Figure \ref{fig:bursting}. The large simulated amplitude modulation of eastern tropical Pacific temperatures $T_2$ provides evidence for a separation of timescales.  \eAR{In the paper by \cite{tim03}}, the simulated El Ni\~no bursting behavior is attributed to a Shilnikov saddle-focus mechanism.
In this paper we seek a more robust mechanism for MMOs, and the analysis relies 
on a global time-scale separation.  Such a global separation can only be seen in 
a dimensionless model, \eJG{because nondimensionalizing the system gives 
insights into the relative sizes of the time-scale parameters that are 
relevant to the MMO analysis.}

\subsection{The Mathematical Model}
Before taking on the scaling, we introduce a change of variables as follows:
\begin{equation*}
	S=T_2-T_1 \qquad T=T_1-T_r \qquad h= h_1+\underbrace{H-z_0}_{K}=h_1+K
\end{equation*} 
With this change of variables, system \eqref{fullT} takes the form
\begin{equation}  
\label{fullS}
\begin{cases}
\displaystyle \frac{dS}{dt} =-\alpha S +\epsilon \mu S^2 +\zeta \mu S \left(S+T+C-C\tanh \left(\frac{h}{h^*}+\frac{bL\mu}{h^*\beta}S\right) \right) \\[.6em]
\displaystyle \frac{dT}{dt} = -\alpha T -\epsilon \mu S^2\\[.6em]
\displaystyle \frac{dh}{dt} =  -r\left(h-K+\frac{bL\mu}{2 \beta} S \right)
\end{cases}
\end{equation}

In system \eqref{fullS}, $S$ is the temperature difference between the eastern and western Pacific surface water, $T$ is the departure of the western Pacific surface ocean temperature from some reference temperature, and $h$ represents the western Pacific thermocline depth anomaly.

\begin{remark}
Despite the change of coordinates, the new system \eqref{fullS} captures the same dynamics as the original system \eqref{fullT} because the two systems are equivalent.  Indeed, they exhibit the same behavior as the dimensionless system \eqref{fast} that we are preparing to introduce.  For example, \cite{kim2011} show that the original system \eqref{fullT} is sensitive to the radiative-convective equilibrium temperature, $T_r$ .  By translating $T_1$ and $T_2$, the value $T_r$ is no longer explicitly present in the new equations $\eqref{fullS}$, however it appears implicitly in the new parameter $C$ {(see Table 2)}.  We make this relationship explicit when we introduce the dimensionless parameter $c$ \eqref{newParam} that corresponds to $C$.  As the analysis will show, the parameter $c$ determines the relative strength of the super El Ni\~{n}o.  Here we see a benefit of the nondimensionalization: it is not $T_r$ alone that is important, but rather the difference $T_r - T_{r_0}$, {which characterizes the maximum possible temperature range in the model}.  Neither $T_r$ nor $T_{r_0}$ {alone} appear in any of the other dimensionless parameters, implying that their effect on the system is entirely dependent on each other.
\end{remark}

\begin{table} [h!]
\begin{center}
\addtolength{\tabcolsep}{1mm}
\renewcommand{\arraystretch}{1.2}
\begin{tabular}{|c|c|c|}
\hline
\textbf{Parameter}  & \textbf{Value} & \textbf{Unit}  \\
\hline
\hline 
T$_0$=S$_0$& 2.8182 & $^\circ$C \\
\hline 
h$_0$ & 62 & m \\
\hline
t$_0$ & 104.9819 & days\\
\hline
$\delta$ & 0.2625 & 1 \\
\hline
$\rho$ & 0.3224  &1  \\
\hline
a & 6.8927 & 1 \\
\hline
k & 0.4032 & 1 \\
\hline
c & 2.3952 & 1 \\
\hline
\end{tabular}
\caption{Rescaled parameters of the system \eAR{from \citep{tim03}}.} \label{rescparameter}
\end{center}
\end{table}

Now seek the dimensionless version of  the system. We define the following dimensionless variables:
\begin{equation*}
	x = \frac{S}{S_0} \qquad y=\frac{T}{T_0} \qquad z= \frac{h}{h_0} \qquad \tau = \frac{t}{t_0}.
\end{equation*}
If we indicate the new time derivative by a prime, i.e. $' = \frac{d}{d \tau}$,  the dimensionless system then takes the form:
\begin{equation}  \label{fast}
\begin{cases}
x' &=\rho \delta (x^2 - ax) + x \left( x+y+c - c \tanh (x+z) \right)\\
y' &= -\rho \delta (ay + x^2) \\
z' &= \delta (k-z- \frac{x}{2}), 
\end{cases}
\end{equation}
where the key parameters can be expressed as 
{\begin{equation} 
\label{newParam}
\begin{array}{llll}
\displaystyle \delta = \frac{r b L }{\zeta h^* \beta} \hspace{1cm} & \displaystyle \rho = \frac{ \epsilon h^* \beta}{r b L  } \hspace{1cm}
& \displaystyle a = \frac{\alpha b L}{\epsilon \beta h^*} \hspace{1cm} & \displaystyle c =\frac{T_r-T_{r0}}{2 S_0}\\[1.5em]
\displaystyle k =\frac{H-z_0}{h^*} & \displaystyle S_0 =T_0=\frac{h^* \beta}{bL\mu} \hspace{1cm}
& h_0 =h^* & \displaystyle t_0 =\frac{ b L}{\beta \zeta h^*}
\end{array}
\end{equation}
}
The zonal advection feedback is represented by the term $\delta x^2$ and the vertical advection and thermocline feedback are captured through $x \left( x+y+c - c \tanh (x+z) \right)$. The term $-\delta  x/2$ is a representation of the recharging mechanism of the thermocline through wind-stress curl-induced Sverdrup transport. Thermal relaxation of sea surface temperatures in this system is expressed in terms of $-\rho \delta ax, -\rho \delta ay$.  The bursting corresponds to the development of extreme El Ni\~no  events with eastern sea surface temperatures approximating those in the western Pacific ($x>-2$), similar to the bursting in equation \eqref{fullS}. The model simulates a rapid termination of these events into a La Ni\~na state ($x<-4$) and a slow recharge of the system through growing small amplitude ENSO oscillations. Upon reaching a threshold state, the system rapidly develops into an extreme El Ni\~no event and the cycle repeats itself.

The system \eqref{fast} has two global time-scale parameters $\varepsilon_1 = \delta$ and $\varepsilon_2 = \rho \delta$.  The parameter $\rho \delta=\varepsilon/\zeta$ relates the efficiency of the advective process to the efficiency of the upwelling term, and $\delta$ is the relative timescale associated with the thermocline adjustment.  Motivated by the bursting behavior seen in Figure \ref{fig:bursting}, we will assume $\delta \ll 1$ so that the system has one fast variable ($x$) and two slow variables ($y,z$).  {Note that the value of the standard model set-up (Table 2) is $\delta=0.26$, but smaller values can be physically justified, if for instance the recharging timescale increases beyond 400 days, or if the upwelling efficiency $\zeta$ increases.}  In the case where both $\delta \ll 1$ and $\rho \ll 1$, then \eqref{fast} will be a three time-scale system with $x$ fast, $z$ slow, $y$ super-slow.  

\begin{figure}[t]
        \centering
        \includegraphics[width=\textwidth]{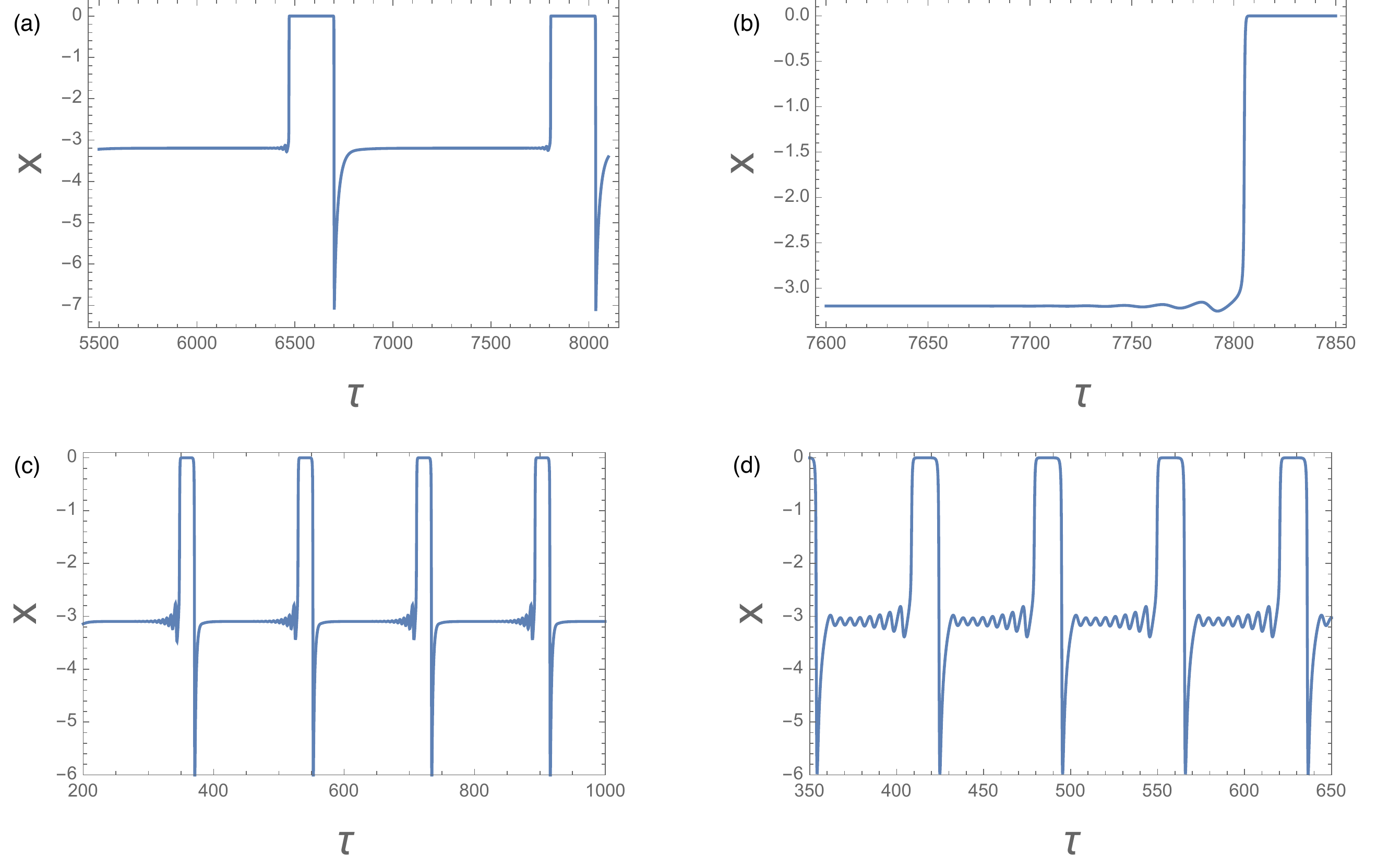}
                \caption{ Examples of MMO patterns in \eqref{fast} when (upper left panel) $\delta=0.01$, $\rho=0.5$, $a = 2.75,$ $c=3.75,$ $k=0.34$, (bottom left panel) $\delta=0.1,$ $\rho=0.35,$  $a = 2.75,$ $c=3.75,$ $k=0.35$, and (bottom right panel) $\delta=0.2,$ $\rho=0.25,$ $a = 2.75,$ $c=3.75,$ $k=0.34$.  Upper right panel: a closer look at the small oscillations in the upper left panel. Notice that the size of the small amplitude oscillations scale with $\delta$, making them difficult to see in the upper left panel.. El Ni\~no events occur when $ x > -2$, whereas strong La Ni\~na events are characterized by $x < -4$. }
        \label{fig:mmoExamples}
\end{figure}

\eCJ{The foundation of the approach we use} for analyzing systems with multiple time scales is called 
{\it Geometric Singular Perturbation Theory (GSPT)} \citep{fenichel79, jones95}.  The combination of GSPT with {\it blow-up techniques} provides
results on global behavior, including complicated patterns such as mixed-mode oscillations \citep{brons, mmoSurvey, ksGSP, ksRO}.  Figure \ref{fig:mmoExamples} depicts mixed mode oscillations observed in equations \eqref{fast}.  Similar bursting behavior has been observed in neurophysiological experiments \citep{amir02, delnegro02, dickson1998, gutfreund1995, khosrovani07}.  Recent mathematical developments have enhanced our ability to explain such phenomena by exploiting a time-scale separation of the underlying model.  The method combines the theory of canards \citep{benoit1983, benoit1981,  szmolwechs, wec07} with a suitable global return mechanism \citep{brons, guckenheimer08, kuehn11, milik1998, aar1, mw05, wec12}, and it is now widely accepted as a robust explanation for \eAR{dynamical} behavior that is qualitatively similar to that observed in Figure \ref{fig:mmoExamples} \citep{brons08, ddw09, mmoSurvey}.

\section{Model Analysis} 
We approach the problem as a singular perturbation problem with respect to the parameter $\delta$, and then relate the dynamics for physically relevant parameters to the geometry of the singular limit.  Thus, we treat the system as a problem with 1 fast variable and 2 slow variables. 

\subsection{GSPT}
Recall that $\tau$ is the {dimensionless} time variable in \eqref{fast}.  We can rescale time by the small parameter $\delta$ by setting $s = \delta \tau$ to obtain the system
\begin{equation}   \label{slow}
\textbf{Slow system}
\begin{cases}
\delta \dot{x} &= \rho \delta (x^2 - ax) + x \left( x+y+c - c \tanh (x+z) \right)\\
 \dot{y} &= -\rho ( ay + x^2 ) \\
\dot{z} &= \left( k-z- \frac{x}{2} \right). 
\end{cases}
\end{equation}
where the dot indicates differentiation with respect to the new variable $s$ (i.e., $\dot{ \ } = d/ds$).  This scaling does not change the trajectories of the system---that is, they have identical phase portraits.  It does 
reparametrize the trajectories, changing the speed at which the curves are traced.  In that sense, the two systems are equivalent as long as $\delta > 0$.  The initial formulation \eqref{fast} is referred to as the `fast system' and \eqref{slow} as the `slow system.'  However, in the limit as $\delta \rightarrow 0$, the systems are different and the case where $\delta = 0$ is referred to as the {\it singular limit}.  

When $\delta=0$, the fast system \eqref{fast} becomes the {\it layer problem}:
\begin{equation}  \label{layer}
\textbf{Layer problem}
\begin{cases}
x' &= x \left( x+y+c - c \tanh (x+z) \right) \\
y' &= 0 \\
z' &= 0. 
\end{cases}
\end{equation}

Physically, this corresponds to a state of constant $T_1$ and $h$, but varying $T_2$. The dynamical evolution of $T_2$ is governed by vertical advective processes in the eastern tropical Pacific.
Notice that the dynamics in the $y$ and $z$ directions are trivial in the layer problem, while the $x$ dynamics are not.  The $y,z$ dynamics are \eJG{described in the \emph{reduced problem}} 
obtained by taking the limit as $\delta\rightarrow 0$ of the slow system \eqref{slow}:
\begin{equation}   \label{reduced}
\textbf{Reduced problem}
\begin{cases}
0 &= x \left( x+y+c - c \tanh (x+z) \right)\\
 \dot{y} &= -\rho ( ay + x^2 ) \\
\dot{z} &= \left( k-z- \frac{x}{2} \right). 
\end{cases}
\end{equation}
This set of equations characterizes a balance between advection and thermocline processes for some value of $(S,T_2,h)$. {Physically this corresponds to a situation in which the} western tropical Pacific SST $T_1$ and thermocline depth $h$ essentially decay towards $T_r$ and $0$, respectively.

\subsection{The Layer Problem: The Critical Manifold and its Stability}
The set
\begin{equation*}
	M = \left\{ x \left[ x+y+c - c \tanh (x+z) \right] = 0  \right\} =\{ F(x,y,z) = 0 \}
\end{equation*}
is important in both components of the singular limit.  In the {\bf layer problem} (equations \eqref{layer}), $M$ is the set of equilibrium points.  As such, the stability of these points for fixed $(y,z)$ values can be computed by the linearization of the  
layer equation. 
\eJG{W}e view $y$ and $z$ as parameters that affect both the location and stability of the equilibria.
\eJG{In} the {\bf reduced problem} (equations \eqref{reduced}), we 
\eJG{observe that $M$ defines the set}
on which the slow ($y,z$) dynamics are defined.  In contrast to a ``standard'' parameter drift approach where $y$ and/or $z$ vary based on an explicit function of time, here $y$ and $z$ vary depending on the fast variable $x$, and hence {\it implicitly} as a function of time.  
\eJG{In the singular limit, the fast system equilibrates to a point in $M$ and then the slow system evolves along $M$.  The importance of $M$ cannot be overstated}---it forms the backbone of complicated dynamical behavior---and it is referred to as the {\it critical manifold}. 

Since the $x'$ equation in \eqref{layer} contains a product, the critical manifold has two components:
\begin{align}
	M_0 &= \{x = 0 \} \\
	M_s &= \left\{ x+y+c - c \tanh (x+z) = 0 \right\} =  \left\{y=  - x - c + c \tanh (x+z) \right\} .
\end{align}
The first condition characterizes an extreme El Ni\~no event ($x=0, T_2=T_1$), depicted e.g. in Figure \ref{fig:mmoExamples}.

GSPT guarantees that for $\delta$ small enough, our intuitive approach will work unless the Jacobian of the layer problem is 0.  Computing the Jacobian, we get
\begin{equation*}
	F_x = x [1 - \sech^2(x+z)] + \left[ x+y+c - c \tanh (x+z) \right].
	\end{equation*}
$M$ is attracting wherever $F_x |_M < 0$ and repelling wherever $F_x |_M > 0.$  First, looking at the stability of $M_0$, we see that 
\begin{equation*}
	F_x |_{M_0} = y+c-c \tanh(z).
\end{equation*}
Therefore $M_0$ is attracting where $y < c \tanh(z) - c$ and repelling when $y > c \tanh(z) - c$.  This condition is crucial in determining the dynamical behavior near an extreme El Ni\~no state ($x=0$).
$F_x = 0$ along the curve $y= c \tanh(z) - c$, and geometrically, this curve represents the intersection of $M_0$ with $M_s$.  Thus, $M_0$ is attracting below $M_s$ and repelling above $M_s$.

\begin{figure}[t!]
        \centering
        \includegraphics[width=\textwidth]{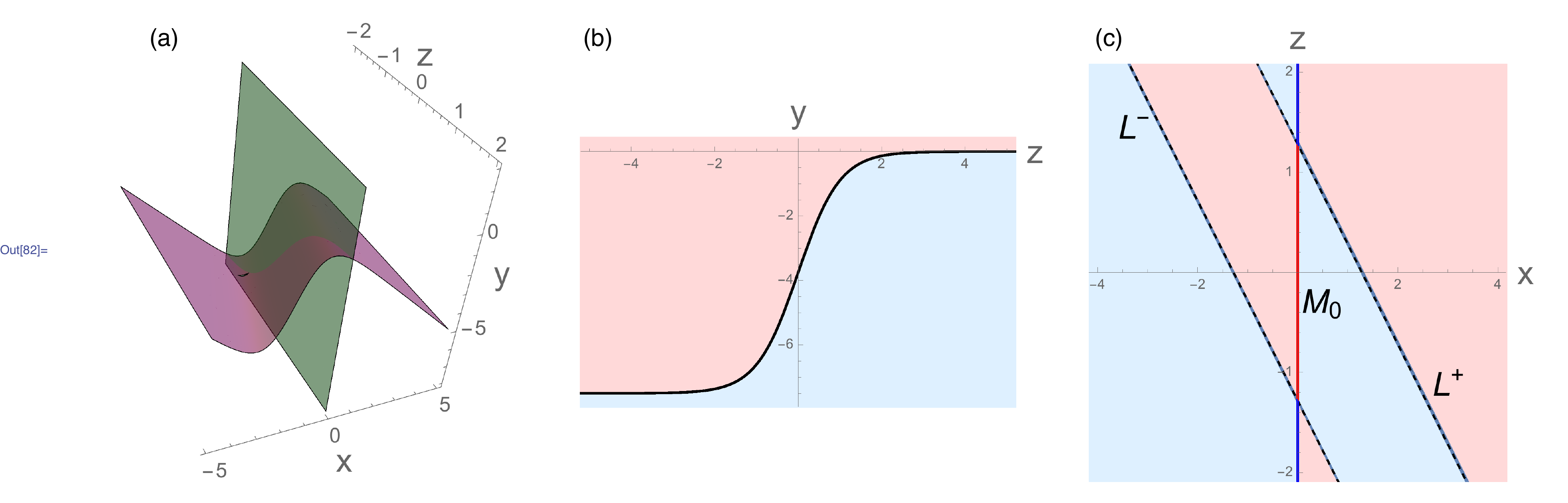}
        \caption{The critical manifold when $a=2.75$, $c=3.75$, $k=0.34$, and $r=0.35$.  Left panel: the manifold in 3D.  Middle panel: the component $M_0$.  Right panel: projection of the component $M_s$ onto the $x,z$-plane. Stable regions are shaded in red, while unstable regions are shaded in blue (middle and right panels only).}
        \label{fig:cmStability}
\end{figure}

Turning to the stability of $M_s,$ we see 
\begin{equation}
	\label{stabilityS}
	F_x |_{M_s} = x[1-c \sech^2(x+z_0)],
\end{equation}
which has three zeroes for fixed $z_0$.  There exist $x_-(z_0) < x_+(z_0)$ such that $c \sech^2(x_\pm+z_0)=1$.  Generically, there are three possible 
\eJG{orders} of these zeros: (1) $x_- < x_+ <0$, (2) $x_- < 0 < x_+$, and (3) $0 < x_- < x_+$.  In 
\eJG{each of these cases}, let $X_1 < X_2 < X_3$ be the roots of $x[1-c 
\sech^2(x+z_0)]$.  Then the left-most branch of $M_s$, i.e. $M_s \cap \{ x < 
X_1\}$, is attracting, with stability alternating at each subsequent root.  Thus $M_s$ is attracting for $x < X_1$ and $X_2 < x < X_3,$ and $M_s$ is repelling 
for $X_1 < x < X_2$ as well as for $x > X_3.$  The lines 
\begin{equation*}
	L^\pm = \{ x=x_\pm(z) \} \cap M_s
\end{equation*}
are called {\it fold curves}.  As seen in Figure \ref{fig:cmStability}, fold 
curves separate the stable and unstable regions of $M_s$.  Along the folds, the standard GSPT approach for normally 
hyperbolic critical manifolds breaks down.  However, it is this degeneracy that 
allows for more complicated interaction between the fast and slow dynamics that 
lie at the heart of the interesting mathematical and {dynamical} behavior 
{of the conceptual ENSO model}.
 
\subsection{The Reduced Problem: Dynamics on the Critical Manifold}
Next, we consider the slow dynamics of the reduced problem \eqref{reduced}.  In 
general, a reduced problem describes a flow of reduced dimension restricted to a 
manifold given by an algebraic condition.  Since $M$ is the union of the 
manifolds $M_0$ and $M_s$, we 
need to solve two separate reduced problems---one on each of the two 
manifolds.  

The reduced problem on $M_0$, described by 
\begin{equation}   \label{reduced0}
	\textbf{Reduced problem on $M_0$}
	\begin{cases}
		0 &= x \\
		 \dot{y} &= -\rho ay \\
		\dot{z} &= k-z. 
	\end{cases}
\end{equation}
is simple because the
global coordinate chart $x = h^0(y,z) = 0$ 
describes the manifold 
$M_0$.  Moreover, the dynamics are already formulated in terms of the coordinate 
variables $y$ and $z$.  Since \eqref{reduced0} is a linear system, it is easy to 
see that $(y,z)=(0,k)$ is the unique globally attracting fixed point on $M_0$ 
(extreme El Ni\~no case).  It is worth noting that, while trajectories on $M_0$ 
will be drawn to this fixed point, it lies in the region where $M_0$ itself is 
unstable.  This observation will play a major role in the global dynamics of 
\eqref{fast} when $\delta>0.$  

\begin{figure}[t]
        \centering
                \includegraphics[width=\textwidth]{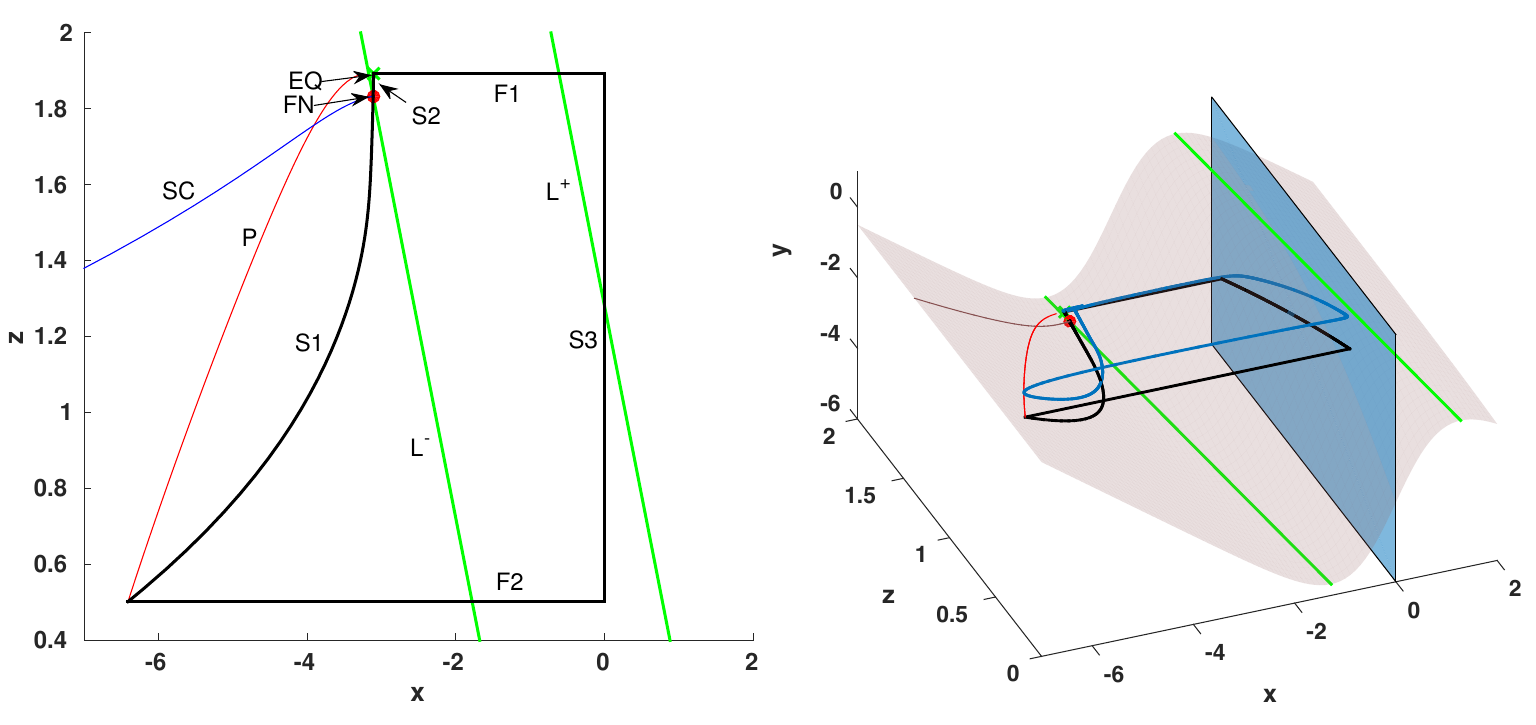}
    	\caption{
      An MMO orbit and approximating singular cycle for parameters
$\delta = 0.1, \rho = 0.5,$ $a = 2.55,$ $c=3.75,$ and $k=0.34.$ Left panel: projection 
of the singular cycle onto the 
$(x,z)$ plane.  The five segments of the singular cycle are 
drawn as black curves and labeled S1, S2, F1, S3 and F2. 
Projections of the fold curves are drawn green and labeled $L^+$ and $L^-$. 
Projection of the segment S3 onto $M_s$ along the fast $x$
direction is drawn as a red curve and labeled P. The folded node 
point is drawn as a large 
red dot and labeled FN. Its strong stable manifold is drawn blue and labeled SC. 
The funnel of FN is the region between SC and $L^-$. The equilibrium point EQ is 
a green x. Right panel: a three dimensional plot
of an MMO orbit (blue) and the approximating 
singular cycle (black). The gray surface is $M_s$ and the light blue surface is $M_0$.}
        \label{fig:enso_sing_cycle}
\end{figure}

Next we turn to the reduced problem on $M_s$, given by
\begin{equation}   \label{reducedS}
\textbf{Reduced problem on $M_s$}
\begin{cases}
0 &=x+y+c - c \tanh (x+z) \\
 \dot{y} &= -\rho ( ay + x^2 ) \\
\dot{z} &= \left( k-z- \frac{x}{2} \right). 
\end{cases}
\end{equation}
In this case we are not as lucky as we were for the problem on $M_0$; the fold lines $L^\pm$ indicate that we will be unable to formulate a global coordinate chart as $x = h^s(y,z)$.  However, using the algebraic condition
\begin{equation*}
	y + x + c - c \tanh (x+z) = 0
\end{equation*} 
we obtain a single coordinate chart, namely
\begin{equation}
\label{alg}
y = h(x,z) = -x - c + c \tanh (x+z),
\end{equation}
where we can study the whole reduced flow on $M_s$.  This is done by differentiating the algebraic condition \eqref{alg} and substituting it for the $\dot{y}$ equation in \eqref{reducedS} to obtain
\begin{align*}
		0 &= \dot{y} -h_x \dot{x} - h_z \dot{z} \\
		\dot{z} &= \left( k-z- \frac{x}{2} \right).
\end{align*}
Substitution and some minor rearranging yields
\begin{equation}
	\label{reduced2}
	\begin{array}{rl}
		[1-c \sech^2 (x+z)] \dot{x}&= \displaystyle \rho (a h(x,z) + x^2) + c \left( k - z - \frac{x}{2} \right) \sech^2(x+z) \\
		\dot{z} &=\displaystyle \left( k-z- \frac{x}{2} \right).
	\end{array}
\end{equation}
Note that the system \eqref{reduced2} is singular along the fold curves $L^\pm$ 
(denoted by green lines in Figure \ref{fig:enso_sing_cycle}) since it is precisely the 
set $\{ 1 - c \sech^2(x+z) = 0 \}$.  

The system \eqref{reduced2} can be desingularized by rescaling the time variable $s$ by a factor of $1 - c \sech^2(x+z)$.  The new system, called the {\it desingularized system}, is
\begin{equation}   \label{desing}
	\textbf{Desingularized problem on $M_s$}
	\begin{cases}
		\dot{x}&= \displaystyle \rho (a h(x,z) + x^2) + c \left( k - z - \frac{x}{2} \right) \sech^2(x+z) \\[1em]
		\dot{z} &=\displaystyle [1-c \sech^2 (x+z)] \left( k-z- \frac{x}{2} \right).
	\end{cases}
\end{equation}

\noindent
The rescaling is quite powerful in that we are now able to define dynamics on 
the entirety of $M_s$, including the folds.  Recalling condition 
\eqref{stabilityS}, we see that the scaling factor changes signs across the fold 
curves $L^\pm$. 
\eJG{The time orientation of trajectories is reversed in the region
where $[1-c \sech^2 (x+z)] < 0$.} 
As indicated in Figure \ref{fig:enso_sing_cycle} we are primarily interested in the region of phase space where $x < 0$.  In fact, $x=0$ is invariant, even for $\delta > 0$ in the full system \eqref{fast}.  Therefore, trajectories starting in the region $x<0$ will remain there for all time.  Physically this means that no El Ni\~no event can grow beyond the warm pool temperature $T_1$.  When $x < 0$, $M_s$ is attracting wherever $[1-c \sech^2 (x+z)] > 0.$  Thus, the rescaling reverses trajectories on the repelling branch of $M_s$ where the time variable is scaled by a negative factor.  

In the desingularized system \eqref{desing}, it is easy to classify three 
types of
\eJG{special points} of the reduced problem:
\begin{itemize}
       \item 
       \eJG{Equilibria} occur where $\dot{x}=0$ and $k - z - x/2 = 
0$.  \eJG{At these points,}
the scaling factor $[1-c \sech^2 (x+z)] \neq 0.$
	\item Regular fold points are fold points that are not equilibria of \eqref{desing}.  That is, $[1-c \sech^2 (x+z)] = 0$ but $\dot{x} \neq 0.$  
	\item Folded singularities are equilibria of \eqref{desing} where $\dot{z} = 0$ due to the rescaling.  That is, folded singularities occur where $\dot{x}=0$ and $[1-c \sech^2 (x+z)] = 0$, but $k - z - x/2 \neq 0$.  The point labeled `FN' in Figure \ref{fig:enso_sing_cycle} is an example of a folded singularity.
\end{itemize}

\eJG{Some folded singularities lead} to the formation of the small amplitude 
oscillations in an MMO.
\eJG{Similar to equilibria}, folded 
singularities are classified by the eigenvalues of the linearization of the flow 
at the equilibrium.  A folded singularity with real, 
\eJG{negative eigenvalues, for 
example, is a stable} folded node, while a folded singularity with complex 
conjugate eigenvalues with positive real part is an unstable folded focus.  The 
`FN' in Figure \ref{fig:enso_sing_cycle} 
\eJG{labels a stable folded node.} 

\subsection{Geometry of MMOs}
\begin{figure}[t]
        \centering
        \includegraphics[width=0.5\textwidth]{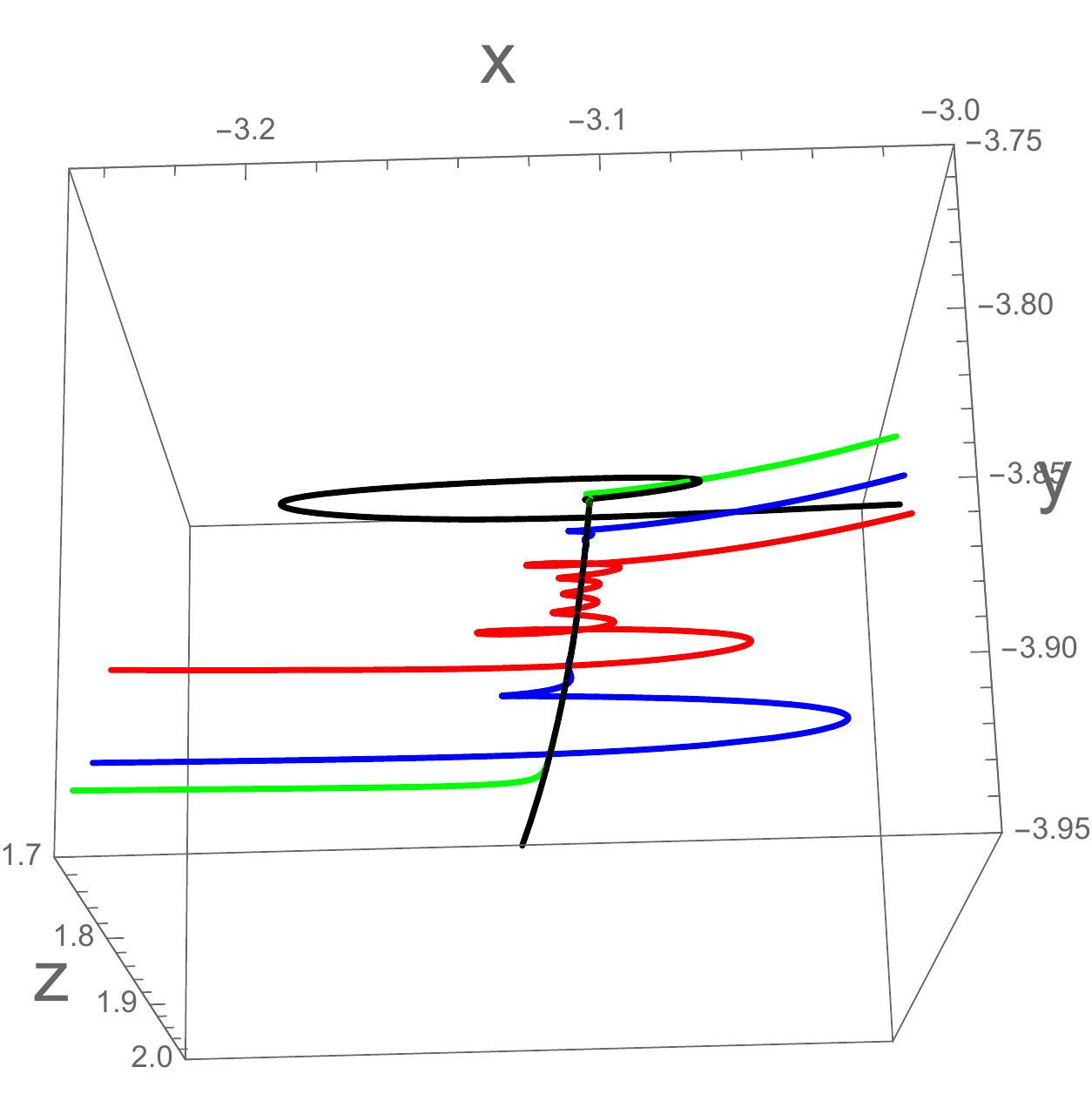}
        \caption{Local rotation (SAOs) of the full system when 
$\delta=0.005$, $\rho = 0.5,$ $a = 2.55,$ $c=3.75,$ and $k=0.34.$  
\eJG{The black trajectory approaches the stable manifold of the saddle-focus  equilibrium along its stable manifold and departs along its unstable manifold. The green, blue, and red
trajectories approach the stable manifold to the equilibrium in the region where the reduced system has a folded node. Their SAOs arise from the twisting of the attracting and repelling slow manifolds in this region. The green trajectory which approaches the stable manifold first has more small oscillations and these oscillations have smaller amplitude than those of the blue and red trajectories.}
}
        \label{fig:twisting}
\end{figure}

\eCJ{\subsubsection{Non-monotone SAOs}}
Folded nodes
can produce MMOs with a suitable global return mechanism \citep{brons, mmoSurvey, guckenheimer08, kw10, kuehn11}.  In a 2D system such as \eqref{desing}, a node is an equilibrium point with two real eigenvalues of the same sign: a weak eigenvalue $\mu_w$ and a strong eigenvalue $\mu_s$ such that $|\mu_w| < |\mu_s|$.  For the folded node denoted `FN' in Figure \ref{fig:enso_sing_cycle}, we have 
\begin{equation}
	\label{evals}\textbf{Eigenvalues at FN}
	\begin{cases}
		\mu_s &\approx -3.48 \\
		\mu_w &\approx -0.13,
	\end{cases}
\end{equation}
 so `FN' is a \eJG{stable} node.
 \eJG{Trajectories approach the node tangent to one of the eigenvectors.  The 
strong stable manifold, corresponding to $\mu_s$, is a single trajectory denoted 
by $\gamma_s$, while all other trajectories approach the folded node in the 
direction of the eigenvector of $\mu_w$. Trajectories that approach the folded 
node from the left without hitting the fold line $L^-$ fill a region called the 
\emph{funnel} of the folded node. The strong stable manifold $\gamma_s$,
denoted `SC' in the left panel of Figure 
\ref{fig:enso_sing_cycle}), 
bounds one side of the funnel. 

The folded node is key to one mechanism for generating SAOs in the ENSO system 
when $\delta >0$. As $\delta$ is increased from 0, the attracting and repelling 
sheets of the critical manifold perturb to invariant slow manifolds that no longer share the same boundaries along the fold lines. The slow manifolds 
twist in a region surrounding the vanished folded node 
point so that they continue to intersect along a finite set of trajectories. The 
number of intersections is related to the eigenvalues of the folded node 
\citep{benoit1991,mw05}. The intersections of the slow manifolds partition the 
attracting slow manifold into \emph{rotational sectors} consisting of 
trajectories that make different numbers of rotations as they pass through the 
folded node region.  \eCJ{A key characteristic of these monotone SAOs is that the number of small amplitude oscillations has an upper bound determined by the geometry of the intersecting slow manifolds.}  Figure \ref{fig:twisting} depicts a few trajectories lying 
in different rotational sectors. Note that as the number of rotations grows, 
their amplitude becomes very small and the oscillations can hardly be seen in the figure 
\citep{mmoSurvey}.} 

\begin{figure}[t]
        \centering
        \includegraphics[width=0.5\textwidth]{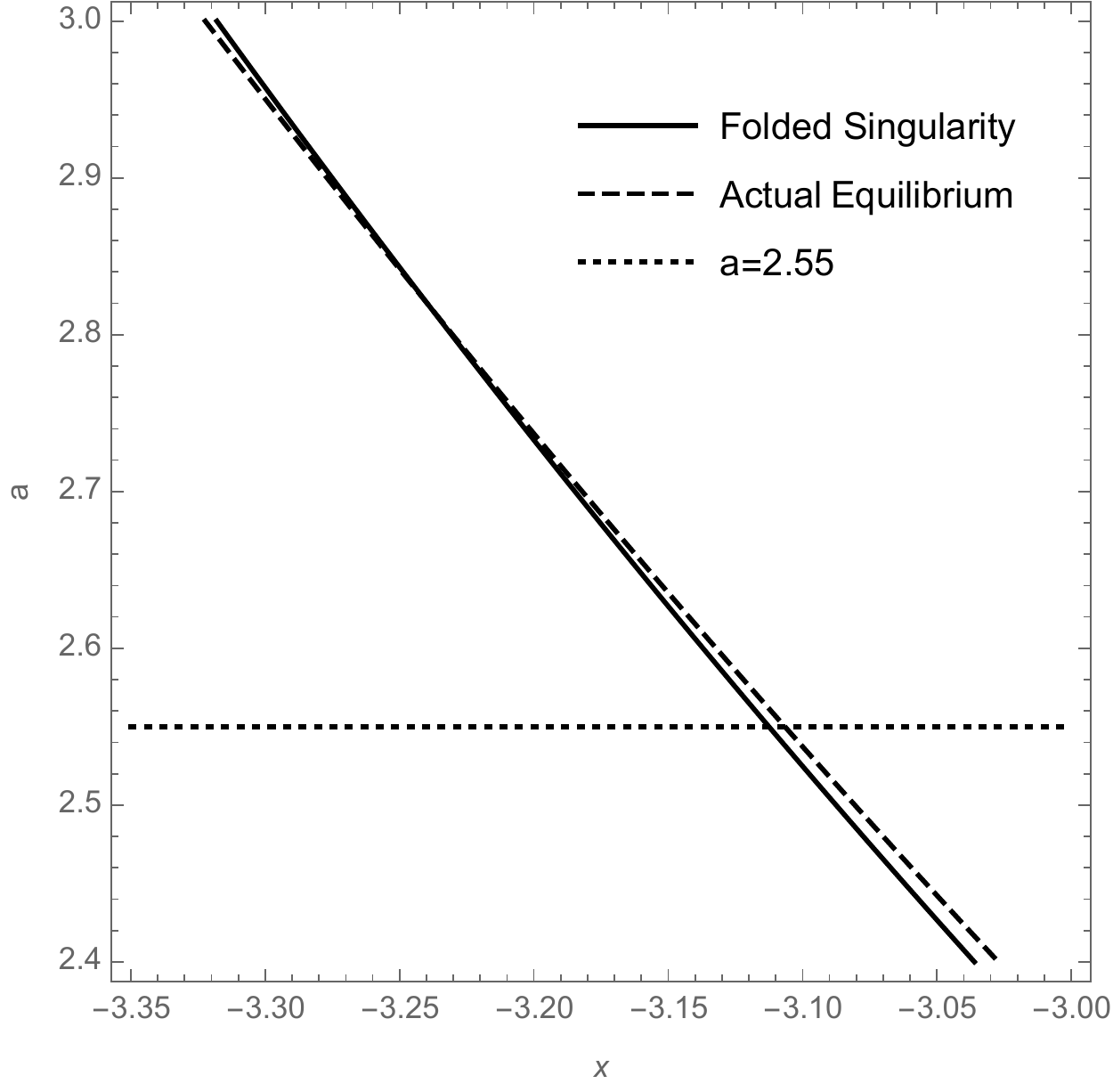}
        \caption{Bifurcation diagram showing how close the system is to a folded 
saddle-node type II when $\rho=0.5,$ $k=0.34$, and $c=3.75$.  A folded 
saddle-node type II occurs when \eqref{fast} \eJG{in its singular 
limit} undergoes a transcritical bifurcation where one equilibrium is an 
ordinary singularity and the other is a folded singularity.}
        \label{fig:fsn2}
\end{figure}

\eCJ{\subsubsection{Monotone SAOs}}
\eJG{The reduced system \eqref{reduced2} has an equilibrium point that lies close to 
the folded node. As parameters are varied in the system, this equilibrium 
crosses the fold curve $L^-$ in a bifurcation that is called a \emph{folded 
saddle-node type II}.  Figure \ref{fig:fsn2} shows the relative location of the equilibrium 
of the singular limit to the folded node as the parameter $a$ varies. 
The folded saddle-node type II bifurcation is located at the 
intersection of these two curves. The equilibrium itself changes from a sink to a saddle as 
it crosses $L^-$. In the full system with $\delta >0$, the equilibrium undergoes 
a \emph{ singular Hopf bifurcation} at parameter values that are $O(\delta)$ 
close to where the equilibrium crosses $L^-$. The adjective singular here refers 
to the fact that the slow and fast time scales both play a role in the 
bifurcation: the imaginary eigenvalues of the equilibrium have magnitude 
$O(1/\sqrt{\delta})$ intermediate between the two time scales. For parameters 
to one side of this Hopf bifurcation, the equilibrium point is a saddle-focus 
with a real negative eigenvalue and a pair of complex eigenvalues with positive 
real parts. Due to the large relative magnitude of the imaginary parts of these 
eigenvalues, trajectories in the unstable manifold of the equilibrium spiral 
away with a large number of oscillations as the amplitude of the oscillations 
grows slowly. These growing oscillations constitute a second mechanism for 
generating SAOs. All trajectories that come close to the equilibrium, leave 
along its unstable manifold and acquire SAOs as they do so. }

\eJG{This behavior connects 
the multiple time scale analysis presented here to the Shilnikov analysis of 
chaotic behavior associated with homoclinic orbits to saddle-focus equilibria in 
three dimensional vector fields \citep{gh83}. Shilnikov bifurcation was a 
mechanism proposed earlier to occur in ENSO 
models of the type studied here \citep{tim03}.}

\eJG{Folded nodes and singular Hopf bifurcation are two mechanisms for producing 
small amplitude oscillations in systems with two slow and one fast variable. For 
these SAOs to be embedded in MMOs, there must be a global return of 
trajectories that leave the SAOs. These global returns may have a singular limit 
as $\delta \to 0$ consisting of a {\it singular cycle} $\Gamma$ that 
concatenates trajectories of the reduced system with trajectories of the layer 
system. Below we find a singular cycle that approximates MMO orbits of the 
ENSO model.}
 
\subsubsection{MMOs in ENSO}

\begin{figure}[t]
        \centering
        \includegraphics[width=\textwidth]{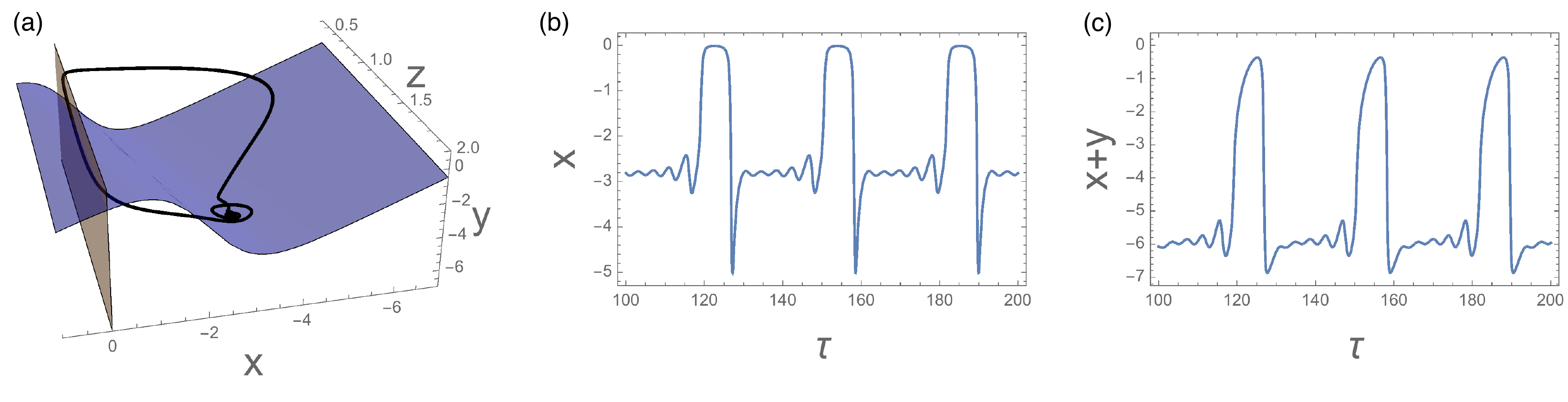}
        \caption{{Attracting MMO when $\delta=0.3$, $\rho = 0.5,$ $a = 2.55,$ 
$c=3.75,$ and $k=0.34.$  Left panel: the MMO orbit in phase space. Middle panel: the time series for $x$ (corresponding to the east/west temperature gradient).  Right panel: the time series for x+y (corresponding to temperature in the eastern Pacific).}}
        \label{fig:MMO3}
\end{figure}

\eJG{Figure \ref{fig:MMO3}a shows a three dimensional phase portrait of an MMO 
in the ENSO model represented by both systems \eqref{fast} and \eqref{slow}. 
Setting $\delta = 0$ produces the singular limit in which trajectories are 
concatenations of solutions to the reduced and  layer equations \eqref{reduced} 
and \eqref{layer}. Here, we construct a singular cycle that approximates this 
MMO orbit and suggests its decomposition into slow and fast segments. 

The singular cycle projected onto the $(x,z)$ plane is shown in the left panel of Figure 
\ref{fig:enso_sing_cycle}. The cycle is composed of five segments. Beginning at the 
lower left corner which represents the {La Ni\~na state ($T_2<T_1$)} of the system following discharge 
of a large El Ni\~no event, the first segment S1 is a trajectory of the reduced 
system on 
$M_s$ that lies in the funnel flowing to the folded node FN. The system 
recharges during this event; as the {western thermocline depth increases}, the temperature in 
the western Pacific $T_1$ increases and the {zonal temperature gradient $(T2-T1)$}
decreases. We anticipate that 
there may be SAOs with decreasing amplitude in this region when $\delta >0$. 
From the folded node, there are a continuum of possible choices for the second 
segment of the singular cycle. (Recall that the desingularization process 
reversed the orientation of the reduced flow on $M_s$ between $L^+$ and $L^-$, 
so that trajectories flow away from the folded node in this strip.) We choose 
the second segment S2 to be the stable manifold of the equilibrium point EQ in the 
desingularized system. This segment is short, and not much appears to happen 
in the phase space except that the trajectory approaches the equilibrium, which 
is a saddle in the reduced system. In the full system, the equilibrium is a 
saddle-focus and SAOs of growing amplitude begin in this region. 

Upon reaching 
the saddle EQ, 
the third segment F1 of the singular cycle is selected to be a jump along a 
trajectory of the layer 
equation to the sheet $M_0$ of the critical manifold. This represents a strong  
El Ni\~no event.  Trajectories do not cross 
$M_0$ since it is invariant in both the full and reduced ENSO models. The fourth 
segment S3 of the singular cycle is then a trajectory of the reduced system for 
$M_0$. Here the thermocline 
and eastern Pacific temperatures evolve on the slow time scale while the 
{zonal temperature gradient} remains small throughout the strong El 
Ni\~no event. Where the singular cycle departs from $M_0$
is not clear: theoretical analysis of what happens 
when two sheets of a slow manifold intersect along a curve is incomplete and a 
subject for future work \citep{thCortez} . The choice we make is to continue 
the trajectory along $M_0$ until the value of $z$ attains the minimum observed on 
the MMO trajectory. From this point, the fifth segment F2 of the singular cycle
is a fast jump back 
to the left hand sheet of $M_s$, discharging the strong El 
Ni\~no event and completing a return to the 
initial point of the cycle. }

\eJG{
The right panel of Figure \ref{fig:enso_sing_cycle} superimposes the MMO computed for parameter values 
$\delta = 0.1, \rho = 0.5,$ $a = 2.55,$ $c=3.75,$ and $k=0.34.$  
on the 
singular cycle displayed in the left panel of Figure \ref{fig:enso_sing_cycle}. The resemblance of the 
two is evident. Using the singular cycle, we can decompose the MMO into slow and 
fast segments. The small amplitude oscillations occur in the region of the reduced
model phase plane close to the folded node FN and the equilibrium point EQ.}  \eCJ{The non-monotone SAOs related to the
folded node mechanism are associated with oscillations in the folded node region, whereas the monotone SAOs due to the singular Hopf mechanism are associated with oscillations that grow along the unstable manifold 
of the equilibrium.}
{The two mechanisms are not mutually exclusive, and determining which is prevalent in ENSO will involve comparison with signals in either real-world observation or output from} {coupled ocean-atmosphere general circulation models.}

\section{Summary and Discussion}
Even though the decadal-scale bursting of strong El Ni\~no events has been the focus of {several} recent studies \cite{kim2011, timmermann02, timmermann03solo}, the underlying dynamics still remain elusive. Our paper provides the first insight into possible connections between MMO dynamics and El Ni\~no bursting. We study the dynamical behavior of a conceptual model of the tropical Pacific climate system, which is based on three nonlinear ordinary differential equations. For realistic physical parameters the system exhibits small amplitude oscillations (SAO), which occasionally burst into large events.  More specifically, the amplitude of the El Ni\~no/La Ni\~na events slowly increases, until reaching a 
\eJG{threshold in} phase space 
where a large El Ni\~no event emerges.\eJG{The event} terminates rapidly as a result of the discharging process (Fig. 1, lower panel). The discharge is so strong, that the system overshoots into a maximum amplitude La Ni\~na event. Its subsequent evolution follows again the growth of SAOs, until the \eJG{threshold} is reached  and the cycle repeats itself.

The underlying equations can be written as a {\it 1 fast, 2 slow system}. The fast variable is associated with the east-west temperature gradient and its rapid adjustment to changes in vertical advection in the eastern tropical Pacific. The slow system captures the adjustment of the western tropical Pacific temperature via thermal relaxation and zonal advection and of the western tropical Pacific thermocline by changes in Sverdrup transport.
To better understand the complex dynamical behavior we apply {\it geometric singular perturbation theory} to the set of conceptual ENSO equations.

In the specific limit of $\delta\rightarrow 0$, which represents either a very slow thermocline discharging process or a very efficient vertical advection, we find that the dynamics of a typical periodic orbit can be spliced together by concatenating trajectories of the decoupled fast and slow dynamics. It is shown that the existence of a folded node, which is characterized by a small El Ni\~no amplitude ($y=-3.095$) plays a key role in organizing the dynamics. This point essentially becomes a {\it point of no return} for the rapid development into a strong El Ni\~no event which follows the fast dynamics. The system can remain on the extreme El Ni\~no manifold ($M_0$) for quite some time (Figs. 2, 5) during which the western tropical Pacific temperature slowly decreases as well as the western tropical thermocline depth. Upon reaching another threshold jump point at western tropical Pacific temperature \eJG{$y \approx -1$}
the system quickly discharges which leads to a rapid transition from $M_0$  onto $M_s$.  From there it moves back slowly towards the folded node to resume its trajectory.

\eCJ{We have shown that a global return mechanism is present in this ENSO Recharge Oscillator model. Associated with this mechanism are two types of MMOs: those with monotone SAOs and those with non-monotone SAOs. The distinction between them can be seen clearly in Figure \ref{fig:twisting}. Both are potential explanations of the interspersion of El Ni\~no bursting with small amplitude oscillations. The key characteristic of the former is the bound on the number of small amplitude oscillations, and for the latter it is the monotonically increasing amplitude of the small oscillations. Whether these can be distinguished in an ENSO signal from observations or large-scale computation is an exciting challenge. }

\section*{Acknowledgments}
A.R., E.W., and C.J. were supported by NSF grant DMS-0940363.  A.T.  was supported by U.S. NSF Grant 10-49219 and U.S. Department of Energy grant DE-SC0005110.  J.G. was supported by U.S. NSF Grant 10-06272.  Additionally, A.R. and E.W. would like to thank UH-Manoa for their hospitality when visiting.

\section*{References}
 \bibliographystyle{ametsoc2014}
      \bibliography{enso_refs}

\end{document}